\newtheorem{thm}{Theorem}

\newtheorem{zam}{Remark}
\newtheorem{pred}{Proposition}

\newtheorem{gip}{Conjecture}

\newtheorem{pro}{Problem}

\documentclass[12pt]{article}
\usepackage{amsfonts}
\usepackage[all]{xy}
\title{On Boris Moishezon's multiple planes
}\author{Vik.S. Kulikov
\thanks{The talk on Workshop on Topology of Algebraic Varieties in honour of Boris Moishezon, Bonn, MPI, Monday July 27, 1998.}}
\date{        }

\begin{document}
\maketitle

\section*{Introduction.  }

Every nonsingular projective surface $S$ over $\mathbb C$ defines three 
underline structures 
$$ vS, \hspace{1cm} dS,\hspace{1cm} tS,$$
where $tS$ is the topological type of $S$, $dS$ is the underline smooth 4-manifold and $vS$ is the deformation type of $S$.

I would like to talk about Boris Moishezon's Program on investigation of smooth structures on projective surfaces and their deformation types.
It has three sources: classical (Italian) algebraic geometry including Picard-Lefschetz theory, braid group theory, topology of smooth manifolds, and it consists of three components. 

The first one coincides with {\it Chisini's Problem}.
Let $S$ be a nonsingular surface in a projective space $\mathbb P^r$ of $\deg S=N$. It is well known that for almost all projections $pr:\mathbb P^r \to \mathbb P^2$ the restrictions $f:S\to \mathbb P^2$ of these projections to $S$ satisfy the following conditions:

$(i)$ $f$ is a finite morphism of $\deg f=\deg S$;

$(ii)$ $f$ is branched along an irreducible curve $B\subset \mathbb P^2$ with ordinary cusps and nodes, as the only singularities; 

$(iii)$  $f^{*}(B)=2R+C$, where $R$ is irreducible and non-singular, and $C$ is reduced; 

$(iv)$  $f_{\mid R}:R\to B$ coincides with the 
normalization of $B$. \newline 
We shall call such $f$ {\it a generic morphism} and its branch curve will be called {\it the discriminant curve}.

Two generic morphisms $(S_1,f_1)$, $(S_2,f_2)$ with the same discriminant curve $B$ are said to be equivalent if there exists an isomorphism $\varphi : S_1 \to S_2$ such that 
$f_1=f_2\circ \varphi $.

The following assertion is known as Chisini's Conjecture. 
\begin{gip}
Let $B$ be the discriminant curve of a generic morphism $f:S\to \mathbb P^2$ of degree $\deg f \geq 5$. Then $f$ is uniquely determined by the pair $(\mathbb P^2,B)$.
\end{gip}

It is easy to see that the similar conjecture for generic morphisms of projective curves to $\mathbb P^1$ is not true. On the other hand one can show that Chisini's Conjecture holds for the discriminant curves of almost all generic morphisms of any projective surface. 

The second part of Moishezon's Program deals with so called {\it braid monodromy technique}. Let $B$ be an algebraic curve in $\mathbb P^2$ of degree $2d$, where $d\in \frac{1}{2}\mathbb N$ (if $B$ is a discriminant curve, then $\deg B$ is even, i.e. $d\in \mathbb N$). The topology of the embedding  $B\subset \mathbb P^2$ is determined by the {\it braid monodromy} of $B$ (see, for example, \cite{Moi1} or \cite{Moi2}), which is described by a factorization of the "full twist" $\Delta _{2d}^2$ in the semi-group $B^{+}_{2d}$ of the braid group $B_{2d}$ of $2d$ 
string braids (in standard generators, $\Delta _{2d}^2=(X_1\cdot ...\cdot X_{2d-1})^{2d}$). If $B$ is a cuspidal curve, then this factorization can be written as follows
\begin{equation}\label{F}
\Delta _{2d}^2=\prod_{i}Q_i^{-1}X_1^{\rho _i}Q_i, \hspace{1cm} \rho _i\in (1,2,3), 
\end{equation} 
where $X_1$ is a positive half-twist in $B_{2d}$.

Let 
\begin{equation} \label{f}
h=g_1\cdot ...\cdot g_r
\end{equation} 
be a factorization in $B^{+}_{2d}$. The transformation which changes two neighboring factors in (\ref{f}) as follows 
$$g_i\cdot g_{i+1} \longmapsto (g_ig_{i+1}g_i^{-1})\cdot g_i, $$
or 
$$g_i\cdot g_{i+1} \longmapsto g_{i+1}(g_{i+1}^{-1}g_ig_{i+1}) $$
is called a {\it Hurwitz move}.

For $z\in B_{2d}$, we denote by 
$$h_z=z^{-1}g_1z\cdot z^{-1}g_2z\cdot ...\cdot z^{-1}g_rz$$
and say that the factorization expression $h_z$ is obtained from (\ref{f}) by simultaneous conjugation by $z$. Two factorizations are called {\it Hurwitz and conjugation equivalent} if one can be obtained from the other by a finite sequence of Hurwitz moves followed by a simultaneous conjugation. For any algebraic curve $B\subset \mathbb P^2$ any two factorizations of the form (\ref{F}) are Hurwitz and conjugation equivalent. We shall say that two factorizations of the form (\ref{F}) belong to the same {\it braid factorization type} if they are Hurwitz and conjugation equivalent. The main problem in this direction is the following one.
\begin{pro}
Does the braid factorization type of the pair $(\mathbb P^2,B)$ uniquely determine the diffeomorphic type of this pair $(\mathbb P^2,B)$, and vice versa?
\end{pro}

Let $S_1$ and $S_2$ be two non-singular projective surfaces, and let 
$\varphi :S_1\to S_2$ be a homeomorphism. The homeomorphism $\varphi $ induces the isomorphism $\varphi ^{*}:H^2(S_2,\mathbb Z)\to H^2(S_1,\mathbb Z)$. Assume that $L_i$, $i=1,2$, is an ample line bundle on $S_i$ such that $f_i:S_i \to \mathbb P^2$ given by three-dimensional linear subsystem of $|L_i|$ is a generic morphism, and let $\varphi ^{*}(L_2)=L_1$. The third part of Moishezon's Program can be formulated as the following problem. 
\begin{pro}
Let $f_i:S_i\to \mathbb P^2$, $i=1,2$, be a generic morphism as above and such that Chisini's Conjecture holds for its discriminant curve $B_i$. Do the diffeomorphic (resp. deformation) types of $S_1$ and $S_2$ coincide if the diffeomorphic (resp. deformation) types of the pairs $(\mathbb P^2,B_1)$ and $(\mathbb P^2,B_2)$ coincide, and vice versa?
\end{pro}

\section{Chisini's Conjecture.}

{\bf 1.1.} Let $B\subset \mathbb P^2$ be an irreducible plane curve  with ordinary cusps and nodes, as the only singularities. Denote by $2d$ the degree of $B$, and let $g$ be the genus of its desingularization, $c= \# \{ \mbox{cusps of}\, B\}$, and $n= \# \{ \mbox{nodes of} \, B\}$.

Let us fix $p\in \mathbb P^2\setminus B$ and denote by $\pi _1=\pi _1(\mathbb P^2\setminus B,\, p)$ the fundamental group of the complement of $B$. Choose any point $x\in B\setminus Sing\, B$ and consider a line $\Pi =\mathbb P^1\subset \mathbb P^2$ intersecting $B$ transversely at $x$. Let $\gamma \subset \Pi$ be a circle of small radius with center at $x$. If we choose an orientation on $\mathbb P^2$, then it defines an orientation on $\gamma $. Let $\Gamma $ be a loop consisting of a path $L$ in $\mathbb P^2\setminus B$ joining the point $p$ with a point $q\in \gamma$, the circuit in positive direction along $\gamma $ beginning and ending at $q$, and  a return to $p$ along the path $L$ in the opposite direction. Such loops $\Gamma$ (and the corresponding elements in $\pi _1$) will be called {\it geometric generators}. It is well-known that $\pi _1$ is generated by geometric generators, and any two geometric generators are conjugated in $\pi _1$ since $B$ is irreducible.

For each singular point $s_i$ of $B$ we choose a small neighborhood $U_i\subset \mathbb P^2$ such that $B\cap U_i$ is defined (in local coordinates in $U_i$) by the equations $y^2=x^3$, if $s_i$ is a cusp, and $y^2=x^2$, if $s_i$ is a node. Let $p_i$ be a point in $U_i\setminus B$. It is well-known that if $s_i$ is a cusp, then $\pi _1(U_i\setminus B,p_i)$ is isomorphic to the braid group $\mbox{Br}_3$ of 3-string braids and is generated by two geometric generators (say $a$ and $b$) satisfying the following relation 
$$aba=bab.$$ 
If $s_i$ is a node, then $\pi _1(U_i\setminus B,p_i)$ is isomorphic to $\mathbb Z\oplus \mathbb Z$ generated by two commuting geometric generators. 

Let us choose smooth paths $\gamma _i$ in $\mathbb P^2\setminus B$ joining $p_i$ and $p$. This choice 
defines homomorphisms $\psi _i:\pi _1(U_i\setminus B,p_i)\to \pi _1$. Denote the image  $\psi _i(\pi _1(U_i\setminus B,p_i))$ by $G_i$ if $s_i$ is a cusp, and $\Gamma _i$ 
if $s_i$ is a node.

A generic morphism of degree $N$ determines a homomorphism $\varphi :\pi _1 \to \mathfrak S_N$, where $\mathfrak S_N$ is the symmetric group. This homomorphism $\varphi $ is determined uniquely up to inner automorphism of $\mathfrak S_N$.
\begin{pred} (\cite{Moi1},\cite{Kul1})\label{prop1}
The set of the non-equivalent generic morphisms of degree $N$ possessing  the same  discriminant curve $B$ is in one to one correspondence with 
the set of the epimorphisms $\varphi :\pi _1(\mathbb P^2\setminus B)\to \mathfrak S_N$ (up to inner automorphisms of $\mathfrak S_N$) satisfying the following conditions:

$(i)$ for a geometric generator $\gamma $ the image $\varphi (\gamma )$ is a transposition in $\mathfrak S_N$; 

$(ii)$ for each cusp $s_i$ the image $\varphi (G_i)$ is  isomorphic to $\mathfrak S_3$ generated by two transpositions;

$(iii)$ for each node $s_i$ the image $\varphi (\Gamma _i)$
is isomorphic to $\mathfrak{S}_2 \times \mathfrak{S}_2$ generated by two commuting transpositions.
\end{pred}
{\bf 1.2.} Moishezon proved the following theorem
\begin{thm} Let an epimorphism $\varphi :\pi _1(\mathbb P^2\setminus B)\to \mathfrak S_N$ satisfy conditions (($i$)-($iii$)) of Proposition \ref{prop1}. If the kernel $K$ of $\varphi $ is a solvable group, then Chisini's Conjecture holds for $B$.
\end{thm} 
In particular, from this theorem, it follows that Chisini's Conjecture holds for the discriminant curves of generic morphisms of $S$ for $S=\mathbb P^2$, $S=\mathbb P^1\times \mathbb P^1$, and for the discriminant curves of generic projections of hypersurfaces in $\mathbb P^3$.
\newline {\bf 1.3.} In \cite{Kul1}, one can find the proof of the following 
\begin{thm} \label{K1}
Let $B$ be the discriminant curve of a generic morphism
$f:S\to \mathbb P^2$ of $\deg f = N$. If 
\begin{equation}
N > \frac{4(3d+g-1)}{2(3d+g-1)-c}. \label{in}
\end{equation}
Then the generic morphism $f$ is uniquely determined by the pair $(\mathbb P^2,B)$ and thus, the Chisini Conjecture holds for $B$.
\end{thm} 

Theorem \ref{K1} shows that if the degree of a generic morphism with given discriminant curve $B$ is sufficiently large, then this generic morphism is unique for $B$. Almost all generic morphisms interesting from algebraic geometric point of view satisfy this condition. More precisely, we have the following theorems (\cite{Kul1}).
\begin{thm}
Let $S$ be a projective non-singular surface, and $L$ be an ample divisor on $S$, $f:S\to \mathbb P^2$ a generic morphism given by a three-dimensional subsystem $\{ E\} \subset |mL|$, $m\in \mathbb Q$, and $B$ its discriminant curve. Then there exists a constant $m_0$ (depending on $L^2, \, (K_S,L), \, K^2_S,\, p_a$) such that $f$ is uniquely determined by the pair $(\mathbb P^2,B)$ if $m\geq m_0$.
\end{thm} 
In particular, we have 
\begin{thm}
Let $S$ be a surface of general type with ample canonical bundle $K_S$, $f:S\to \mathbb P^2$ a generic morphism given by a three-dimensional linear subsystem of $|E|$, where $E\equiv mK_S$, $m\in \mathbb N$ ($\equiv $ means numerical equivalence). Then $f$ is uniquely determined by the pair $(\mathbb P^2,B)$.
\end{thm} 
{\bf 1.4.} The last theorem can be generalized to the case of generic morphisms $f:X\to \mathbb P^2$, where $X$ is a canonical model of a surface $S$ of general type ($X$ is the surface with Du Val singularities) (\cite{Kul4}).
\begin{thm} (V.S. Kulikov and Vik.S. Kulikov)
Let $S_1$ and $S_2$ be two minimal models of surfaces of general type with $K_{S_1}^2=K_{S_2}^2$ and $\chi (S_1)=\chi (S_2)$, and $X_1$ and $X_2$ their canonical models. Let $B$ be the canonical discriminant curve of two $m$-canonical generic morphisms $f_1:X_1\to \mathbb P^2$ and  $f_2:X_2\to \mathbb P^2$, that is, $f_i$, $i=1,2,$ is given by three-dimensional linear subsystems of $|mK_{X_i}|$. If $m\geq 4$, then $f_1$ and $f_2$ are equivalent.
\end{thm}

\section{Braid factorization and smooth equivalence.}
{\bf 2.1.} One can show that the braid factorization of an algebraic 
plane curve $B$ uniquely determines the diffeomorphic type of the pair 
$(\mathbb P^2,B)$. More precisely, one can prove
\begin{thm} (Vik.S.Kulikov and M.Teicher)
Let $B_1, B_2 \subset \mathbb P^2$ be two projective plane curves with the same braid factorization type. Then there exists a diffeomorphism 
$\varphi : \mathbb P^2 \to \mathbb P^2$ such that $\varphi (B_1)=B_2$.
\end{thm}
The inverse question remains open. 
\begin{pro} \label{pr}
Let $B_1, B_2 \subset \mathbb P^2$ be two projective plane curves such that the pairs $(\mathbb P^2,B_1)$ and $(\mathbb P^2,B_2)$ are diffeomorphic. Do $B_1$ and $B_2$ belong to the same braid factorization type?
\end{pro}

I would like to mention here the following very important problem.
\begin{pro} 
Let $\Delta^2_{2d}={\cal{E}}_1$ and $\Delta^2_{2d}={\cal{E}}_2$ be two braid factorizations. Does there exist a finite algorithm to recognize whether these two braid factorizations belong to the same braid factorization type or not?
\end{pro}

\section{Smooth types of surfaces and smooth types \\ of pairs $(\mathbb P^2,B)$.}
{\bf 3.1.} 
The set of plane curves of degree $2d$ is naturally parameterized by the points in $\mathbb P^{d(2d+3)}$. The subset of plane irreducible curves of degree $2d$ and genus $g$ with $c$ ordinary cusps and some nodes, as the only singularities, corresponds to a quasi-projective subvariety ${\cal{M}}(2d,g,c)\subset P^{d(2d+3)}$ (\cite{Wah}). One can show that if two non-singular points of the same irreducible component of ${\cal{M}}_{red}(2d,g,c)$ correspond to curves $B_1$ and $B_2$, then the pairs $(\mathbb P^2,B_1)$ and $(\mathbb P^2,B_2)$ are diffeomorphic. In particular, in this case the fundamental groups $\pi _1(\mathbb P^2\setminus B_1)$ and $\pi _1(\mathbb P^2\setminus B_2)$ are isomorphic. Moreover, in this case $B_1$ and $B_2$ have the same braid factorization type.

The following Proposition is a simple consequence of Proposition \ref{prop1} and some local properties of generic morphisms.

\begin{pred}
Let $(\mathbb P^2,B_1)$ and $(\mathbb P^2,B_2)$ be two diffeomorphic (resp. homeomorphic) pairs. If $B_1$ is the discriminant curve of a generic morphism $(S_1,f_1)$, then 
$B_2$ is also the discriminant curve of some generic morphism $(S_2,f_2)$. Moreover, if  $(S_1,f_1)$ is unique, i.e. Chisini's Conjecture holds for $B_1$, then the same is true for $(S_2,f_2)$ and  $S_1$ and $S_2$ are diffeomorphic (resp. homeomorphic).
\end{pred}
{\bf 3.2.} Let $B$ be the discriminant curve of some generic morphism $f:S\to \mathbb P^2$ of $\deg f = N$. Suppose Chisini's Conjecture holds for $B$. Let $f^{*}(B)=2R+C$ be the preimage of $B$ and suppose that for $\pi (S\setminus C)$ there exists a unique (up to inner automorphism) epimorphism $\varphi :\pi (S\setminus C)\to \mathfrak{S}_{N-1}$. In this case, we say that Chisini's Conjecture holds {\it twice} for $B$. 
If it is the case, one can prove the following 
\begin{thm} (\cite{Kul3}) 
Let $B_i$, $i=1,2,$ be the discriminant curve of a generic morphism $f_i:S_i\to \mathbb P^2$. Put $f_i^{*}(B_i)=2R_i+C_i$. Suppose  Chisini's Conjecture holds twice for $B_i$. Then there exists a diffeomorphism (resp. a homeomorphism) of pairs $\varphi : (\mathbb P^2,B_1) \to  (\mathbb P^2,B_2)$ if and only if 
there exists a diffeomorphism  (resp. a homeomorphism) $\Phi :S_1 \to S_2$ such that $\Phi(R_1)=R_2$, $\Phi(C_1)=C_2$ and such that the following diagram is commutative
$$
\xymatrix{
S_1 \ar[d]_{f_1} \ar[r]^{\Phi } & S_2 \ar[d]^{f_2} \\
\mathbb P^2 \ar[r]_{\varphi } & \mathbb P^2 . }
$$
\end{thm}
{\bf 3.3.} For $S\subset \mathbb P^r$, a projection $f:S\to \mathbb P^2$ is defined by a point in Grassmannian $\mbox{Gr}_{r+1,r-2}$ (the base locus of the projection). It is well known that the set of generic projections is in one to one correspondence with some Zariski's open subset $U_S$ 
of $\mbox{Gr}_{r+1,r-2}$. A continuous variation of a point in $U_S$ gives rise to a continuous family of generic projections of $S$, whose branch curves belong to the same continuous family of plane cuspidal curves. Therefore the discriminant curves of two generic projections of $S\subset \mathbb P^r$ belong to the same irreducible component of ${\cal{M}}(2d,g,c)$. Moreover, it is easy to see that the discriminant curves of two generic projections of $S\subset \mathbb P^r$ have the same braid factorization type, because they belong to the same irreducible component of ${\cal{M}}(2d,g,c)$ .
    
In particular, if two surfaces $S_1$ and $S_2$ of general type with the same $K^2_S=k$ and $p_a=\chi ({\cal{O}}_{S_i}) = p$ are embedded by the $m$th canonical class into the same projective space $\mathbb P^r$ and belong to the same 
irreducible component of coarse moduli space ${\cal{M}}_S(k,p)$ of surfaces with given invariants (\cite{Gie}), then there exist generic projections $f_1$ of $S_1$ and $f_2$ of $S_2$ belonging to the same continuous family of generic projections. Therefore, discriminant curves (we will call them {\it $m$-canonical discriminant curves}) of two such generic projections of 
$S_1$ and $S_2$, belonging to the same 
irreducible component of a moduli space ${\cal{M}}_S(k,p)$, belong to the same irreducible component of ${\cal{M}}(2d,g,c)$ (cf. \cite{Wah}).
By Theorem 2 and by Propositions 5 and 6 in \cite{Kul1}, for a surface of general type with ample canonical class the triple of integers $(m,k,p)$ is uniquely determined by the invariants $(d,g,c)$ of $m$th canonical discriminant curve, and vice versa. Thus, we have a natural mapping $ir_{k,p,m}$ (resp. $var_{k,p,m}$) from the set of irreducible (resp. connected) components of ${\cal{M}}_S(k,p)$ to the set of irreducible (resp. connected) components of ${\cal{M}}(2d,g,c)$.  Hence by Proposition 1 and Theorems 4 and 5, we have 
\begin{thm} (\cite{Kul4}) Assume that for all surfaces $S$ of general type with the same $K^2_S=k$ and $p_a=\chi ({\cal{O}}_{S_i})=p$ there exists a three dimensional linear subsystem of $|mK_S|$ which gives a generic morphism $f:S\to \mathbb P^2$. Then
\newline ($i$) $ir_{k,p,m}$ is injective for any $m\in \mathbb N$; \newline ($ii$) $var_{k,p,m}$ is injective if $m\geq 4$. 
\end{thm} 
{\bf 3.4.} 
From Theorems 4, 6  and Proposition 2 it follows
\begin{thm} (Vik.S. Kulikov and M. Teicher) Let $B_i$, $i=1,2,$ be  
an $m$-canonical discriminant curve of a generic morphism $f_{m,i}:S_i\to \mathbb P^2$. If $B_1$ and $B_2$ have the same braid decomposition type, then $S_1$ and $S_2$ are diffeomorphic.
\end{thm}

Applying Theorems 4 and 5 we obtain 
\begin{thm}
The braid factorization type of an $m$-canonical discriminant curve $B_m$ is an invariant of the deformation type of the corresponding surface $S$ of general type.
\end{thm}
\begin{zam}
A negative solution of Problem \ref{pr} implies a negative solution of Diff-Def Problem.
\end{zam}
\begin{zam}
In \cite{Man}, Manetti announced that Diff-Def Problem has a negative solution.
\end{zam}

{\bf 3.5.}
One can prove the following theorem using the arguments similar to those in \cite{Kul1}.
\begin{thm} (\cite{Kul3}) Chisini's Conjecture holds twice for the $m$-canonical discriminant curves.
\end{thm}
Applying this theorem and theorem 7, we have 
\begin{thm} \label{dif} (\cite{Kul3}) 
Let $B_i$, $i=1,2,$ be an $m$-canonical discriminant curve of some generic morphism $f_{m,i}:S_{i}\to \mathbb P^2$. Put $f_i^{*}(B_i)=2R_i+C_i$.  Then the pairs $(\mathbb P^2,B_1)$ and $(\mathbb P^2,B_2)$ are diffeomorphic (resp. homeomorphic) if and only if there exists a diffeomorphism  (resp. a homeomorphism) $\Phi :S_1 \to S_2$ such that $\Phi(R_1)=R_2$, $\Phi(C_1)=C_2$ and such that the following diagram is commutative
$$
\xymatrix{
S_1 \ar[d]_{f_1} \ar[r]^{\Phi } & S_2 \ar[d]^{f_2} \\
\mathbb P^2 \ar[r]_{\varphi } & \mathbb P^2 , }
$$ 
where $\varphi $ is a diffeomorphism  (resp. a homeomorphism) of pairs $(\mathbb P^2,B_i)$ .
\end{thm}
{\bf 3.6.}
In \cite{Cat1} and \cite{Cat2}, Catanese investigated smooth simple bidouble coverings $\varphi : S\to Q=\mathbb P^1\times \mathbb P^1$ of type $(a,b),(m,n)$; thus, $\varphi $ is a finite $(\mathbb Z/2)^2$ Galois covering branched along two generic curves of respective bidegrees $(2a,2b)$, $(2m,2n)$. He proved that for each integer $k$ there exists at least one $k$-tuple $S_1,...,\, S_k$ of bidouble coverings of $Q$ of respective types $(a_i,b_i),(m_i,n_i)$ satisfying the following conditions:  

$(i)$ $S_i$ and $S_j$ are homeomorphic for each $1\leq i,j\leq k $.

$(ii)$ $r(S_i)\neq r(S_j)$ for $i\neq j$ (and therefore $S_i$ and $S_j$ are not diffeomorphic), where $r(S_i)= \max \{ s\in \mathbb N\, \, |\, \, (1/s)K_{S_i}\in H^2(S_i,\mathbb Z)\} $ is the index of $S_i$. \newline We shall call a $k$-tuple of surfaces of general type,  satisfying conditions ($i$) and ($ii$), {\it a Catanese $k$-tuple}.

As a corollary of Theorem \ref{dif} we obtain 
\begin{thm} (\cite{Kul2}) 
Let $B_{m,i}$, $i=1,2,$ be $m$-canonical discriminant curves of generic morphisms $f_{m,i}:S_{i}\to \mathbb P^2$, where $(S_1,S_2)$ is a Catanese pair. Then the pairs $(\mathbb P^2,B_{m,1})$ and $(\mathbb P^2,B_{m,2})$ are not homeomorphic. In particular, $B_{m,1}$ and $B_{m,2}$ belong to different braid factorization types.
\end{thm} 

Theorem 13 shows that the braid factorization types of $m$-canonical discriminant curves of Catanese's surfaces distinguish diffeomorphic types of these surfaces.

\section{Moishezon 4-manifolds.}
{\bf 4.1.} In \cite{Moi1}, Moishezon remarked that for any cuspidal factorization of $\Delta^2_{2d}=\cal{E}$ and a projection $pr:\mathbb C^2\to \mathbb C^1$, it is possible to construct a cuspidal curve $\overline B$ and a topological embedding $i:\overline B\to \mathbb P^2$ which is not complex-analytic, but behaves as a complex-analytic one with respect to the rational map $\pi: \mathbb P^2\to \mathbb P^1$, so that the braid monodromy for $i(\overline B)$ with respect to $\pi $ will be well defined and will be represented by $\Delta^2_{2d}=\cal{E}$. We shall call such curves $i(\overline B)$ {\it semi-algebraic}. 
Assume that there exists an epimorphism $\varphi :\pi _1(\mathbb P^2\setminus \overline B)\to \mathfrak S_N$ satisfying the conditions described in Proposition 1. Then one can construct a smooth 4-manifold $S$ and 
a $\cal{C}^{\infty}$-map $f:S\to \mathbb P^2$ which outside of 
$i(\overline B)$ is a non-ramified covering of degree $N$ with monodromy homomorphism equal to $\varphi$, and over a neighborhood of $i(\overline B)$ locally behaves like a generic morphism branched at $i(\overline B)$ and with local monodromy induced by $\varphi$. We call such $S$ {\it the Moishezon 4-manifold}.
\begin{gip} The class of Moishezon 4-manifolds coincides with the class of symplectic 4-manifolds.
\end{gip}
{\bf 4.2.} {\bf Problem}. {\it To realize Boris Moishezon's Program for the Moishezon 4-manifolds.}

Steklov Mathematical Institute

{\rm victor$@$olya.ips.ras.ru }

\end{document}